\UseRawInputEncoding
\documentclass[12pt]{article}

\usepackage[all]{xy}
\usepackage{bm}
\usepackage{amsfonts}
\usepackage{mathrsfs}
\usepackage{amsfonts}
\usepackage{amsfonts}
\usepackage{amsfonts}
\usepackage{amsfonts}
\usepackage{amsfonts}
\usepackage{amsfonts}
\usepackage{amsfonts}
\usepackage{amsfonts}
\usepackage{amsfonts}
\usepackage{amsfonts}
\usepackage{amsfonts}
\usepackage{amsfonts}
\usepackage{amsfonts}
\usepackage{amsfonts}
\usepackage{amscd}
\usepackage{bbm}
\usepackage [dvips]{graphics}
\usepackage[centertags]{amsmath}
\usepackage{amsfonts}
\usepackage{amssymb}
\usepackage{amsmath}
\allowdisplaybreaks
\usepackage{xypic}
\usepackage{amsthm}
\usepackage{amsmath}
\usepackage{newlfont}
\usepackage{stmaryrd}
\usepackage[]{titletoc}  
\usepackage{mathrsfs}
\usepackage{stmaryrd}
\usepackage [dvips]{graphics}
\usepackage [none]{hyphenat}
\usepackage[centertags]{amsmath}
\usepackage{amsfonts}
\usepackage{amssymb}
\usepackage{amsthm}

\usepackage{newlfont}
\usepackage{hyperref}

\newcommand{\bx}{$\hfill\Box$}
\newlength{\defbaselineskip}

\newcommand{\setlinespacing}[1]%
           {\setlength{\baselineskip}{#1 \defbaselineskip}}

\theoremstyle{plain}
\newtheorem{thm}{Theorem}[section]

\newtheorem{cor}[thm]{Corollary}
\newtheorem{lem}[thm]{Lemma}
\newtheorem{exa}[thm]{Example}

\newtheorem{prop}[thm]{Proposition}

\newtheorem{rem}[thm]{Remark}
\newtheorem{Def}[thm]{Definition}

\newcommand{\cc}{\mathbb{C}}

\newcommand{\la}{\lambda}

\makeatletter\@addtoreset{equation}{section} \makeatother
\begin{document}

\textwidth=13.5cm
  \textheight=23cm
  \hoffset=-1cm

  \baselineskip=17pt

\title { Pl$\ddot{\mathrm{u}}$ckerians twisted with linear forms and Druzkowski map }
\author{Li Chen}
\date{}
\maketitle {\textbf{}
\def\zz{\textbf{z}}
\def\pkk{Pl$\ddot{{\mathrm{u}}}$cker }
\def\ww{\textbf{w}}
\def\zzz{\textbf{z}_0}
\def\bg{\mbox{\boldmath$\gamma$}}
\def\ba{\mbox{\boldmath$a$}}
\def\bb{\mbox{\boldmath$b$}}
\def\bet{\mbox{\boldmath$\beta$}}
\def\M{\mathcal{M}}
\def\H{\mathcal{H}}
\def\ff{\footnotesize}
\def\nn{\normalsize}
\allowdisplaybreaks[1]

 \textbf{Abstract:} Let  $C$ be an $n\times n$ matrix  such that  $ 3\leq n\leq 2\, \mathrm{rank}\,C$ and $T$ be the  cubic linear  map     with respect to  $C$.  We introduce and apply an algebraic construction (called Pl$\ddot{\mathrm{u}}$cker  polynomials) to show that  if     the Jacobian determinant of $T$ is constant, then the row matroid of $C$ is non-uniform, and the condition $3\leq n\leq 2\, \mathrm{rank}\,C$ can not be relaxed. This is  one  step forward from the fundamental $J(T)$=const$\Longrightarrow\det C=0$ constraint,  and the key ingredient in the proof is the interplay between the standard \pkk relation and an unexpected    linear rigidity phenomenon obtained via the  Pl$\ddot{\mathrm{u}}$cker  polynomials.   We also exhibit independent interests of these polynomials as    a variation of the classical Pl$\ddot{\mathrm{u}}$cker-Grassmann construction.



\section{Introduction and main results}
\subsection{ Constant Jacobian  and Matroid structure  }

 Let $T$ be a polynomial map and $J(T)$ be its Jacobian determinant.  As  Alp$\ddot{\mathrm{o}}$ge's recent   counterexample disproves the Jacobian Conjecture in dimension greater than 2, the problem that to what extend can  $J(T)$=const constrains  the global structure of $T$,  besides invertibility,   remains   illusive and widely open.   In this paper we focus on the following cubic linear maps, which encodes all non-linear parts of a polynomial map into a single matrix.
 \begin{Def}


Let $C=(c_{ij})_{1\leq i,j\leq n}$  be an $n\times n$   matrix over $\mathbb{C}$. The map  $ T: \footnotesize (z_1,\cdots,z_n)\mapsto (w_1,\cdots,w_n)$  \normalsize   defined by
\begin{equation}\label{jutii}  w_i=z_i+f_i^3, \quad  f_i=c_{i1}z_1+c_{i2}z_2+\cdots +c_{in}z_n,1\leq i\leq n  \end{equation}\normalsize  is called a \emph{cubic linear map} (or Dru$\dot{z}$kowski map). \end{Def}
Cubic linear maps are  well-known  to be  the stabilizing model  for general constant-Jacobian polynomial maps by the fundamental works of  Bass-Connell-Wright  and Dru$\dot{z}$kowski \cite{BCW,LD1}
 hence has been extensively studied \cite{Ch,debondtyan2013,LD2,hubb,TM,Es,DW}.

The most basis Jacobian determinant forced     constraint on (\ref{jutii}) is well-known to be  $$J(T)=\mathrm{const}\Rightarrow\mathrm{rank}\,C<n. $$ We  will go further from this to  consider    how $J(T)$=const  constrains the matroid structure of $C$ via the following notion of  general and degenerate row position.



 \begin{Def}\label{nizhen} A    rank $r$ matrix  $C$     is said to be in general  row position if any $r$ rows of $C$   are linearly independent. Otherwise,  $C$     is said to be in degenerate   row position.

 \end{Def}
In standard combinatorial language, general(resp. degenerate) row position means that the row matroid of $C$ is     uniform (resp. non-uniform).
 Our main result as follows  states that if $\mathrm{rank }\,C$ is sufficiently large, then $J(T)$=const forces non-uniform matroid structure on $C$.
\begin{thm}[\textbf{Degeneracy Theorem}]\label{badai}Let  $C$ be a complex $n\times n$ matrix  such that $3\leq n\leq 2\,\mathrm{rank} \, C$. If the Jacobian determinant       of the cubic linear map (\ref{jutii})      is a   constant, then $C$ is  in degenerate row position.

\end{thm}
The following example shows that    the condition $3\leq n\leq 2\,\mathrm{rank} \, C$   can not be dropped.
 \begin{exa}To see $n\leq 2\,\mathrm{rank}\,C$ is sharp, let   $C=\scriptsize\begin{pmatrix}
1 & -1 & 0 \\
1 & -1 & 0 \\
1 & -1 & 0
\end{pmatrix} $. Then $C$ is of rank 1 hence $n=2\, \mathrm{rank}\,C+1$. One easily checks that the map (\ref{jutii})   has Jacobian determinant 1, but    $C$ is in general row position since each row is nonzero. Similarly, let $C=\scriptsize\begin{pmatrix}
1 & -1   \\
1 & -1
\end{pmatrix} $, then  $J(T)$=1, $n=2\,\mathrm{rank }\,C$ and $C$ is in general position, which shows  $n\geq 3$ is sharp.\end{exa}


We outline  the proof of Theorem \ref{badai}.  For any $k$-index $ (i_1,\cdots,i_k)$, let  $  m_{i_1i_1\cdots i_k}   $  denote  the $k\times k$ principal minor of $C$ with respect to $i_1,\cdots,i_k$.  It is well-known that the Jacobian determinant of  (\ref{jutii}) equals  the following sum over   $1\leq k\leq n$ where each $k$  corresponds to a homogeneous component of degree $2k$
  \begin{equation}\label{33333}J(T) = 1 + \sum_{k=1}^{n} 3^k \sum_{{1\leq i_1<i_2<\cdots <i_k\leq n}} m_{i_1i_1\cdots i_k}f_{i_1}^2f_{i_2}^2\cdots f_{i_k}^2. \end{equation}
 We will introduce a new algebraic construction called \pkk polynomials (Definition \ref{pksum}) to prove  a rigidity phenomenon (Theorem \ref{jiefangcheng}) on linear relations of $\{f_{i_1}^2f_{i_2}^2\cdots f_{i_k}^2\}$, which   applies to extract the precise  values(up to a global scalar multiple) of   $\{m_{i_1i_1\cdots i_k}\}$ from   $J(T)$=const, and the proof will be completed by combining these values  with  the classical three-term \pkk relation  via full rank   decomposition of $C$.

  \subsection{\pkk polynomials and  linear rigidity  }

Relations among the minors of a matrix are   an important topic in many areas of algebra \cite{Fulton1991,Holtz2007,Lin2009,Oeding2011,Procesi1996,Weyman2003}. The classical Pl$\ddot{\mathrm{u}}$cker theory  gives quadratic relations on $l\times l$ minors of an $n\times l$ matrix, and after removing overlapping indices, the relations can be reduced to those on $2l\times l$ matrices. Here   we    introduce    a  new Pl$\ddot{ {\mathrm{u}}}$cker-type expression called \pkk polynomials for a $2l\times l$ matrix,  which involves both its $l\times l$  minors and row linear forms.
 \begin{Def}\label{pksum}Let $l$ be a positive integer,     $M=\{a_{ij}\}_{1\leq i\leq 2l, 1\leq j\leq l}$ be a $2l\times l$ matrix of rank $l$  and $g_i=a_{i1}z_1+a_{i2}z_2+\cdots+a_{il}z_l$,  $1\leq i\leq 2l$.   The  sum
  \begin{equation}\label{pmii}  P_M:=\sum_{I=(i_1,i_2,\cdots, i_l)}(-1)^{|I|+1}D_ID^3_{I^c}g_{i_1}^2g_{i_2}^2\cdots g_{i_l}^2  \quad \end{equation} over all  $\binom{2l}{l}$ indices $\{I=(i_1,i_2\cdots, i_l)|1\leq i_1<\cdots <i_l\leq 2l\}$ is called the   {Pl$\ddot{ {u}}$cker   polynomial} with respect to $M$. Here
      $|I|=i_1+\cdots+i_l$, $I^c$ is     the complementary index of  $I$ in $(1,2,\cdots,2l), $  and $D_I$ is  the determinant of the      matrix  consisting of   $i_1,\cdots,i_l$ rows of $M$.



 \end{Def}
Non-triviality of Definition \ref{pksum} starts from $l=2$(see Example \ref{2lizi}). The involvement of $D_I$ and $D_{I^c}$ in (\ref{pmii}) is natural as they also  appear in the classical $G(l,2l)$ \pkk relation, while the $\{g_{i_1}^2g_{i_2}^2\cdots g_{i_l}^2\}$ part  is directly motivated by the constant Jacobian condition (\ref{33333}).  The    $ {|I|+1}$ power on $-1$ and cubic power on $D_{I^c}$ guarantee  that $P_M$ is homogeneous and anti-symmetric respect to rows of $M$(see Lemma \ref{hanglieshi} below).

 The suprising thing  we  find is that if $M$ is in general row position, then $P_M=0$ is exactly the   criterion for   linear dependence of $\{g_{i_1}^2g_{i_2}^2\cdots g_{i_l}^2\}$, and their   nontrivial linear relation  is uniquely determined by coefficients of $P_M$.



 \begin{thm}[\textbf{Rigidity Theorem}]\label{jiefangcheng}Let $M$ be a    $2l\times l$ matrix of  rank $l$   which is   in general row position,  and $\{g_i\}_{i=1}^{2l}$ be row linear forms of $M$. The followings are equivalent

   $(i)$ $P_M=0$.

 $(ii)$ The linear equation \begin{equation}\label{zhuyaode}\sum_{1\leq i_1<i_2<\cdots <i_l\leq 2l}x_{i_1 i_2\cdots  i_l}g_{i_1}^2g_{i_2}^2\cdots g_{i_l}^2=0     \end{equation} with  $\binom{2l}{l}$ scalar unknowns $\{x_{i_1i_2\cdots i_l}\}$  admits nontrivial solutions.

  In this case, the solution space is one dimensional spanned by the {Pl$\ddot{ {u}}$cker}  polynomial coefficients $\{(-1)^{|I|+1}\lambda D_ID^3_{I^c}\}$.
   \end{thm}
In matroid language, the last statement implies  that the product set $\{g_{i_1}^2g_{i_2}^2\cdots g_{i_l}^2\}$ is  a circuit in the sense that no proper subset of it can be   dependent. The following  example on $l=2$ helps to visualize Definition \ref{pksum} and Theorem \ref{badai}.
%
%
 \begin{exa}\label{2lizi}Write $M$   as    $\scriptsize\left(
\begin{array}{cccc}
 a_1 & b_1 & c_1 & d_1 \\
 a_2 & b_2 & c_2 & d_2 \\
\end{array}
\right)^T$\normalsize, then  $P_M$   equals the LHS  of the following algebraic identity(we will prove this identity  as part of  Example \ref{36})
     \begin{equation}\label{22nnb}
\text{\footnotesize  $\displaystyle
\begin{aligned}
  &\quad(a_1b_2 - a_2b_1) (c_1d_2 - c_2d_1)^3 (a_1z_1 + a_2z_2)^2 (b_1z_1 + b_2z_2)^2 + (c_1d_2 - c_2d_1) (a_1b_2 - a_2b_1)^3 (c_1z_1 + c_2z_2)^2 (d_1z_1 + d_2z_2)^2 \\
&- (a_1c_2 - a_2c_1) (b_1d_2 - b_2d_1)^3 (a_1z_1 + a_2z_2)^2 (c_1z_1 + c_2z_2)^2 - (b_1d_2 - b_2d_1) (a_1c_2 - a_2c_1)^3 (b_1z_1 + b_2z_2)^2 (d_1z_1 + d_2z_2)^2 \\
&+ (a_1d_2 - a_2d_1) (b_1c_2 - b_2c_1)^3 (a_1z_1 + a_2z_2)^2 (d_1z_1 + d_2z_2)^2 + (b_1c_2 - b_2c_1) (a_1d_2 - a_2d_1)^3 (b_1z_1 + b_2z_2)^2 (c_1z_1 + c_2z_2)^2\\&=0
\end{aligned}$}
\end{equation}  hence all      $4\times 2$ Pl$\ddot{ {u}}$cker  polynomials are  identically zero.  On the other hand,  the 6 polynomials  $\{(a_1z_1 + a_2z_2)^2 (b_1z_1 + b_2z_2)^2 ,\cdots\}$ in (\ref{22nnb}) can be spanned by 5 monomials $\{z_1^4,z_1^3z_2,z_1^2z_2^2,z_1z_2^3,z_2^4\}$ hence  automatically  dependent. This verifies the $(i)\Leftrightarrow(ii)$ part of  Theorem \ref{badai}  on $l=2$ where both statements hold unconditionally.
  \end{exa}
Besides Theorem \ref{jiefangcheng} as an important application,  we will  initiate a   study on the structure of  $P_M$   itself
which  is of independent interest as a variation  of the classical Pl$\ddot{\mathrm{u}}$cker-Grassmann theory. In fact, even in case $l=2$, the mechanism underlying identity  (\ref{22nnb}) is already quite different from that of the  classical $G(2,4)$ \pkk relation
$$\left(a_1 b_2-a_2 b_1\right)\left(c_1 d_2-c_2 d_1\right)-\left(a_1 c_2-a_2 c_1\right) \left(b_1 d_2-b_2 d_1\right)+
 \left(a_1 d_2-a_2 d_1\right) \left(b_1 c_2-b_2 c_1\right)=0$$ though they are both of the $[ab][cd]-[ac][bd]+[ad][bc]=0$ form. When $l\geq 3$,
  $P_M$ is no longer identically zero and we will exhibit    an inductive      structure  upgrade of $P_M$  as   $l$ goes up, standing in contrast to standard \pkk theory where matrix size is of much less  conceptual   importance.
  We will prove Theorem \ref{jiefangcheng} and Theorem \ref{badai} in Section 2, and the general structure of $P_M$ will be studied in Section 3.



\section{Proof of the main results}
We first apply  Theorem \ref{jiefangcheng} to  prove Theorem \ref{badai} in Section 2.1,  then prove Theorem \ref{jiefangcheng} in Section 2.2.
\subsection{Proof of Theorem \ref{badai}}
Decomposing (\ref{33333}) into homogeneous components shows  that $J(T)$=const if and only if for every $1\leq k\leq n$,
  \begin{equation}\label{zhenfanrena} \sum\nolimits_{1\leq i_1<i_2<\cdots <i_k\leq n}m_{i_1i_2\cdots i_k}f_{i_1}^2f_{i_2}^2\cdots f_{i_k}^2=0.   \end{equation}
 In light of the following   elementary  proposition   relating row matroid structure of $C$ to principal minors of $C$ of size $\mathrm{rank}\,C\times \mathrm{rank}\,C$,  we will prove Theorem \ref{badai} by exploring  (\ref{zhenfanrena}) specified to $k=\mathrm{rank}\,C$.

 \begin{prop}\label{mzf}Let $C$ be an $n\times n$ matrix with rank $r<n$ and $C=M_{n\times r}\widetilde{M}_{r\times n}$ be a full rank decomposition of $C$.

$(i)$ $C$ is in general row position if and only if $M$ is in general row position.

$(ii)$ Let $I=(i_1,\cdots,i_r)$ be an $r$-index and $m_I$ be the principal minor of $C$ with respect to $I$, then $m_I=D_I\widetilde{D}_I$,  where $D_I$ is the $r\times r$ minor of $M$ consisting of rows from $I$ and $\widetilde{D}_I$ is the $r\times r$ minor of $\widetilde{M}$ consisting of columns from $I$.

$(iii)$ If all $r\times r$ principal minors of $C$ vanish, $C$ is in degenerate row position.

\end{prop}

\begin{proof}$(i) $ For any $I=(i_1,\cdots,i_r)$, denote the submatrices of $C$ and $M$ using rows from $I$ by $C_I$ and $M_I$, then   $C_I=M_I\widetilde{M}$. Since $\widetilde{M}$ is of rank $r$,   $C_I$ is of rank $r$ if and only if the square matrix $M_I$ is non-singular.

$(ii)$ Obvious.

$(iii)$ If $C$ is in general row position, then $M$ is in general row position by $(i)$ hence  $D_I\neq 0$ for all $r$-indices $I$. This  combined with $(ii)$ and $m_I=0$  forces all $\widetilde{D}_I=0$, which   is a contradiction   since $\widetilde{M}$ is of full rank.

\end{proof}
The following proof of Theorem \ref{badai}  splits into $n<2\mathrm{rank}\,C$ and $n=2\mathrm{rank}\,C$. The case $n=2\mathrm{rank}\,C$ turns out to   carry  the main difficulty where Theorem \ref{jiefangcheng} will be needed.
\subsubsection{The cases $n<2\,\mathrm{rank}\,C$}

%
We now prove all $n<2\,\mathrm{rank}\,C$ cases of Theorem \ref{badai}.

\begin{proof}Let $C$, $\{f_i\}_{i=1}^n$ be as in (\ref{jutii})  and write $r:=\mathrm{rank}\,C$.  Suppose Theorem \ref{badai} is false, that is,  $C$ is in general row position. We show that  the $\binom{n}{r}$ polynomials $$\{f_{i_1}^2f_{i_2}^2\cdots f_{i_r}^2|1\leq i_1<i_2<\cdots <i_r\leq n\}$$ are linearly independent. Once this is proved,   (\ref{zhenfanrena}) specialized to $k=r$ implies vanishing of all $r\times r$ principal minors  of $C$, and Proposition \ref{mzf} $(iii)$ in turn  forces degenerate row position of $C$, which is a  contradiction and thus completes the proof.

Now suppose that
\begin{equation}\label{xwg}\sum_{I } x_If_{I}^2=0\end{equation} for scalars $\{x_I\}$ (here $f_I$ stands for $f_{i_1}f_{i_2}\cdots f_{i_r}$). We show that for any  fixed $r$-index $I_0$, $x_{I_0}=0$.
Let $S=\{1,2,\cdots,n\}\setminus I_0$ which is a  nonempty set of cardinality $n-r$, and set $$Z=\{\zz=(z_1,\cdots,z_n)|f_k(\zz)=0,\forall k\in S\}.$$ Then for any $r$-index $I$ different from $I_0$, $f_I|_Z=0$ and thus (\ref{xwg}) restricted to $Z$ gives $$ x_{I_0}f_{I_{0}}^2|_Z=0$$ and it remains to show that  $f_{I_{0}}^2|_Z$ is not identically zero, which will force  the desired $x_{I_0}=0.$
Note that  $\{f_k|k\in S\}$ consists of  linear forms,  the ideal $$I(Z)=\{p\in\cc[z_1,\cdots,z_n]|p(\zz)=0,\forall \zz\in Z\} $$  is   prime  hence it suffices to show that for any $i\in I_0$, $f_i|_Z$ is not identically zero, and this follows from $n<2r$. In fact, if $f_i|_Z$ is identically zero, then $f_i$ lies in the linear space spanned by $\{f_k|k\in S\}$ so $\{f_i\}\cup \{f_k|k\in S\}$ is a linearly dependent set of $\{f_1,\cdots,f_n\}$ with cardinality $n-r+1$. Since $n<2r$, this cardinality is at most $r$, which is impossible since  $C$ is  in general row position. The proof is completed. \end{proof}
One sees that above argument   fails exactly when  $n$ hits $2\,\mathrm{rank}\,C$, where general row position of $C$ does not imply linear independence of products of the squared linear forms,  hence the principle minors can indeed be nonzero. We will prove  the $n=2\,\mathrm{rank}\,C$  case of Theorem \ref{badai}   by a conceptually different mechanism.


\subsubsection{The case $n=2\,\mathrm{rank}\,C$}
We will still assume $C$ is in general position and  draw a contradiction to complete the proof. Though in case  $n=2\,\mathrm{rank}\,C$, this no longer  implies vanishing of the principle minors, we can  invoke Theorem \ref{jiefangcheng} which gives explicit values of the principal minors, and  the desired contradiction is obtained  by showing that these values violate  the classical three-term \pkk relation.

\begin{lem}\label{3xiangmu}(Three-term Pl$\ddot{ {u}}$cker relation) Let $l\geq 2$ and $M$ be a  $2l\times l$ matrix. Let   $S \subset \{5, 6, \dots, 2l\}$ be a fixed subset of cardinality  $ l - 2$ (when $l=2$, $S$ is taken to be the empty set). Then it holds that
\begin{equation}\label{eq:plucker_1234}
D_{12S} D_{ 34S} - D_{ 13S} D_{ 24S} + D_{ 14S} D_{ 23S} = 0.
\end{equation}
Here $D_{ijS}$ denotes the determinant of the $l\times l$ submatrix of $M$ consisting rows whose index set is $S \cup \{i, j\}$.

\end{lem}
\begin{lem}\label{lifanghe}
Let $a,b,c$ be complex numbers. If  $a-b+c=0$ and $a^3-b^3+c^3=0$, then $abc=0$.

\end{lem}

\begin{proof}Using $b=a+c$, we have
$$0=a^3-b^3+c^3=a^3-(a+c)^3+c^3=-3ac(a+c)=-3abc.$$
\end{proof}
\begin{rem} The conceptual point in Lemma \ref{lifanghe} is that it fails when $a,b,c$ lie  in a field of   characteristic 3,  where $3abc=0$ does not imply $abc=0$.  By   (\ref{33333}), characteristic 3 is exactly the case  where    $J(T)$=1  always holds,  hence nothing on $C$ can be said from it.

\end{rem}
Now we  prove the  $n=2\,\mathrm{rank}\,C$ case of  Theorem \ref{badai}.

\begin{proof}Write $n=2l$ then   $\mathrm{rank}\,C=l$. Moreover,  $l\geq 2$ since $n\geq 3$ as assumed in Theorem \ref{badai}. Let $C$, $\{f_i\}_{i=1}^{n}$ be as in (\ref{jutii}) and  
\begin{equation}\label{guanxii}
C=M_{2l\times l}\widetilde{M}_{l\times 2l}. \end{equation} be a full rank decomposition of $C$. Let  $D_I$ be the $l\times l$ minor of $M$ consisting of its rows from $I$ and  $\widetilde{D}_I$ be the $l\times l$ minor of $\widetilde{M}$ consisting of its columns from $I$.

Assume that the conclusion of Theorem \ref{badai} is false, that is, $C$ is in general row position,  then by Proposition \ref{mzf} $(i)$, $M$ is in general row position, that is,
\begin{equation}\label{bulinga}D_I\neq 0   \end{equation} for every $l$-index $I$.  We draw a contradiction from this to complete  the proof.

Let $\{g_i\}_{i=1}^{2l}$  be row linear forms of $M$. Then by (\ref{guanxii}), $\{g_i\}_{i=1}^{2l}$ and $\{f_i\}_{i=1}^{2l}$ are related by
$$f_i=g_i\circ \widetilde{M}$$ which combined with
(\ref{zhenfanrena}) specified to $k=l$ gives
 \begin{equation}\label{zhenfanrenaa} \sum\nolimits_{1\leq i_1<i_2<\cdots <i_l\leq 2l}m_{i_1i_2\cdots i_l}(g_{i_1}\circ \widetilde{M})^2(g_{i_2}\circ \widetilde{M})^2\cdots (g_{i_l}\circ \widetilde{M})^2=0.  \end{equation}
 Since $\widetilde{M}$ is of full rank, above equation gives
 \begin{equation}\label{zhenfanrenaaa} \sum\nolimits_{1\leq i_1<i_2<\cdots <i_l\leq 2l}m_{i_1i_2\cdots i_l} g_{i_1}^2 g_{i_2}^2\cdots  g_{i_l}^2=0   \end{equation} which means that $\{m_I\}$ is a   solution to equation (\ref{zhuyaode}). Since we assume $C$ is in general row position, this solution must be nontrivial by Proposition \ref{mzf} $(iii)$.


By Theorem \ref{jiefangcheng}, there exists a nonzero scalar $\lambda$ such that
$$m_I=(-1)^{|I|+1}\lambda D_ID^3_{I^c}$$ which combined with $m_I=D_I\widetilde{D}_I$(by Proposition \ref{mzf} $(ii)$) and $D_I\neq 0$ gives
\begin{equation}\label{jibenshangle}
\widetilde{D}_I=(-1)^{|I|+1}\lambda D^3_{I^c}.
\end{equation}
Now  take a  fixed subset $$S \subset \{5, 6, \dots, 2l\}$$   of cardinality  $ l - 2$ and let $$T=\{5,6,\cdots,2l\}\setminus S.$$ Here  both    $S$ and $T$ are  taken to be the empty set if $l=2$. Specializing   (\ref{jibenshangle}) to the six $l$ indices $12S,13S,14S,23S,24S,34S$ gives
\begin{equation}\label{123456}
\begin{split}
\widetilde{D}_{12S} = (-1)^{|12S|+1}\lambda D^3_{34T}, \quad
\widetilde{D}_{13S} = (-1)^{|13S|+1}\lambda D^3_{24T}, \quad
\widetilde{D}_{14S} = (-1)^{|14S|+1}\lambda D^3_{23T}, \\
\widetilde{D}_{23S} = (-1)^{|23S|+1}\lambda D^3_{14T}, \quad
\widetilde{D}_{24S} = (-1)^{|24S|+1}\lambda D^3_{13T}, \quad
\widetilde{D}_{34S} = (-1)^{|34S|+1}\lambda D^3_{12T}.
\end{split}
\end{equation}

By Lemma \ref{3xiangmu}, it holds that

\begin{equation}\label{dag24}
D_{12T}D_{34T}-D_{13T}D_{24T}+D_{14T}D_{23T}=0
\end{equation}
\begin{equation}\label{dag34}
\widetilde{D}_{12S}\widetilde{D}_{34S}-\widetilde{D}_{13S}\widetilde{D}_{24S}+\widetilde{D}_{14S}\widetilde{D}_{23S}=0
\end{equation}

Inserting (\ref{123456}) into (\ref{dag34}) and discarding  the nonzero $\lambda$  gives 

\begin{equation}\label{wuao}
D^3_{12T} D^3_{34T} - D^3_{13T} D^3_{24T} + D^3_{14T} D^3_{23T} = 0
\end{equation}
where    all powers of $-1$ cancel (for instance $\widetilde{D}_{12S}\widetilde{D}_{34S}=(-1)^{|12S|+1+|34S|+1}\lambda D^3_{34T}\lambda D^3_{12T}=\lambda^2 D^3_{34T}D^3_{12T}$ since the power of $-1$ is even).

 Combining  (\ref{dag24}), (\ref{wuao}) and Lemma \ref{lifanghe} gives
 $$D_{12T}D_{34T}D_{13T}D_{24T}D_{14T}D_{23T}=0,$$ which  contradicts (\ref{bulinga}). The proof is complete.
\end{proof}

 \subsection{Proof of Theorem \ref{jiefangcheng}}

\begin{proof}
$(i)\Rightarrow(ii)$. This is obvious,  as $P_M$   is now  a specific zero linear sum of
$\{g_{i_1}^2g_{i_2}^2\cdots g_{i_l}^2\}$, which is nontrivial by the assumption $D_I\neq 0$ for all $I$.

$(ii)\Rightarrow(i)$. It suffices to prove the last statement of the theorem. That is, given      a nontrivial zero linear sum
\begin{equation}\label{0he}\sum_I x_Ig^2_I=0,  \end{equation} there exists a constant $\la\neq0$ such that \begin{equation}\label{shima}x_I=\la(-1)^{|I|+1}D_ID^3_{I^c}\end{equation} for every $l$-index $I$.
Write $c_I:= (-1)^{|I|+1}D_ID^3_{I^c}$ for short, then (\ref{shima}) follows if we can show that
\begin{equation}\label{zhenshi}\frac{x_I}{c_I}=\frac{x_J}{c_J},  \forall I\neq J.
\end{equation}
In what follows,  we show that (\ref{zhenshi}) holds
if    $I$ and $J$ have  exactly    $l-1$ common entries, and this will complete the proof since  any two $l$-indices  can be connected by finitely many indices where every adjacent pair  have $l-1$ common entries.

 By symmetry, we  prove  (\ref{zhenshi}) assuming  $$I=(1,2,\cdots,l) $$ and $J$ lies in the following index set

\begin{equation}\label{iij}\{I_{ij}:=(1,2,\cdots,\hat{i},\cdots,l,l+j), 1\leq i,j\leq l\} \end{equation} consisting  of all $l^2$  indices   sharing  $l-1$ entries with $I$, that is,

\begin{equation}\label{zhenshia}\frac{x_I}{c_I}=\frac{x_{I_{ij}}}{c_{I_{ij}}}.
\end{equation}

\bigskip
 We   first prove that (\ref{zhenshia}) holds when $M$ is of the standard block form $\scriptsize\left(
                              \begin{array}{c}
                                A \\
                                I \\
                              \end{array}
                            \right). $
 Write $A=\{a_{ij}\}$, then $\det(A)\neq 0$ as the top $l\times l$ minor. Moreover, each $a_{ij}$ equals   some $l\times l$ minor consisting of one row from $A$ and $l-1$ rows from the lower identity block up to a sign, hence $a_{ij}\neq 0$ for all $1\leq i,j\leq l$.

 Our strategy is to fix \(j\in\{1,\ldots,l\}\) and   trace the coefficients of the following $l$ monomials
$$
z_j^{2l} \quad \mathrm{and }\quad
 \{z_j^{2l-1}z_k|1\leq k\leq l, k\neq j\}
$$
in (\ref{0he}). It is easy to see that
only    \(l+1\) terms
$$
x_Ig_I^2
 \quad
\mathrm{and}
\quad
 \{x_{I_{ij}}g^2_{I_{ij}}|1\leq i\leq l\}
$$
in (\ref{0he})  can contribute to these coefficients (keep in mind that $I=(1,2,\cdots,l)$), which will produce  a linear system consisting of $l$ equations with $l+1$ variables $
\{x_I,x_{I_{1j}},x_{I_{2j}},\cdots,x_{I_{lj}}\},
$ and we will  prove (\ref{zhenshia}) by solving this linear system.



\bigskip

For the coefficient of  $z_j^{2l}$,  the term
$
 x_Ig_I^2
$
contributes to
$
x_I\prod_{i=1}^l a_{ij}^2
$
and each term
$x_{I_{ij}}g^2_{I_{ij}}(=x_{I_{ij}}z_j^2\prod\limits_{\substack{1 \le r \le l \\ r \ne i}} g_r^2)$ contributes to  $x_{I_{ij}}\prod\limits_{\substack{1 \le r \le l \\ r \ne i}}a_{rj}^2$.
Now (\ref{0he}) gives
\begin{equation}\label{liangl}
x_I\prod_{i=1}^l a_{ij}^2+\sum_{i=1}^l(x_{I_{ij}}\prod\limits_{\substack{1 \le r \le l \\ r \ne i}}a_{rj}^2)=0.
\end{equation}
For notational simplicity, we write
 \begin{equation}\label{daititi}P_j:=\prod_{i=1}^l a_{ij}^2  \quad \mathrm{and} \quad  u_i:=x_{I_{ij}}\prod\limits_{\substack{1 \le r \le l \\ r \ne i}}a_{rj}^2,1\leq i\leq l,\end{equation}
 then (\ref{liangl}) becomes
\begin{equation}\label{liangll}
x_IP_j+\sum_{i=1}^lu_i=0.
\end{equation}

The monomial  $z_j^{2l-1}z_k,k\neq j$ appears in $g_I^2$ when one $g_i^2$ contributes to $z_jz_k$ via the mixed term
$
2a_{ij}a_{ik}z_jz_k,
$  while other $g_r^2,r\neq i$ contribute  to $z_j^2$ via $a_{rj}^2z_j^2$, hence it has coefficient $2x_I
\sum_{i=1}^l(
a_{ij}a_{ik}
\prod\limits_{\substack{1 \le r \le l \\ r \ne i}} a_{rj}^2) $.
Similarly,   $z_j^{2l-1}z_k$ has coefficient
$
2x_{I_{ij}}
\sum\limits_{\substack{1 \le r \le l \\ r \ne i}}
(a_{rj}a_{rk}
\prod\limits_{\substack{1\leq s\leq l\\s\neq r,i}}a_{sj}^2
)$ in $x_{I_{ij}}g^2_{I_{ij}}$ for each $1\leq i\leq l$.  Now (\ref{0he}) gives

\begin{equation} x_I
\sum_{i=1}^l(
a_{ij}a_{ik}
\prod\limits_{\substack{1 \le r \le l \\ r \ne i}} a_{rj}^2)+\sum_{i=1}^l( x_{I_{ij}}
\sum\limits_{\substack{1 \le r \le l \\ r \ne i}}
a_{rj}a_{rk}
\prod\limits_{\substack{1\leq s\leq l\\s\neq r,i}}a_{sj}^2)=0, k\neq j.\end{equation}
%
%


Invoking above notations $P_j$ and $u_i$, this becomes
\begin{equation}\label{liang111}
x_IP_j
\sum_{i=1}^l
\frac{a_{ik}}{a_{ij}}
+
\sum_{i=1}^l
u_i
\sum\limits_{\substack{1 \le r \le l \\ r \ne i}}
\frac{a_{rk}}{a_{rj}}
=0,k\neq j.
\end{equation}
Using (\ref{liangll}), equation (\ref{liang111}) further simplifies into (after replacing $\sum\limits_{\substack{1 \le r \le l \\ r \ne i}}
\frac{a_{rk}}{a_{rj}}$ by  $
(\sum_{r=1}^l\frac{a_{rk}}{a_{rj}})-\frac{a_{ik}}{a_{ij}}$ in (\ref{liang111}))

\begin{equation}\label{liang112}
\sum_{i=1}^l
\frac{a_{ik}}{a_{ij}}u_i=0,k\neq j.
\end{equation}
Viewing  $u_1,u_2,\cdots,u_l$ as unknowns, then combining    (\ref{liang112}) specialized to all   $$k=1,2,\cdots,j-1,j+1,\cdots,l$$ together with (\ref{liangll})
gives  a linear system consisting of  $l$ equations
$$
\begin{pmatrix}
\frac{a_{11}}{a_{1j}} & \frac{a_{21}}{a_{2j}} & \cdots & \frac{a_{l1}}{a_{lj}} \\
\vdots & \vdots & \ddots & \vdots \\
\frac{a_{1,j-1}}{a_{1j}} & \frac{a_{2,j-1}}{a_{2j}} & \cdots & \frac{a_{l,j-1}}{a_{lj}} \\
1 & 1 & \cdots & 1 \\
\frac{a_{1,j+1}}{a_{1j}} & \frac{a_{2,j+1}}{a_{2j}} & \cdots & \frac{a_{l,j+1}}{a_{lj}} \\
\vdots & \vdots & \ddots & \vdots \\
\frac{a_{1l}}{a_{1j}} & \frac{a_{2l}}{a_{2j}} & \cdots & \frac{a_{ll}}{a_{lj}}
\end{pmatrix}
\begin{pmatrix}
u_1 \\
u_2 \\
\vdots \\
\vdots \\
\vdots \\
u_l
\end{pmatrix}
=
\begin{pmatrix}
0 \\
\vdots \\
0 \\
-x_I P_j \\
0 \\
\vdots \\
0
\end{pmatrix}.$$
The determinant of the coefficient matrix of this linear system  equals $ \frac{\det(A)}{\prod_{i=1}^{l} a_{ij}}$ which is nonzero.
Applying Cramer's rule we obtain the unique solution
\begin{equation}\label{cramersolution}u_i = x_I P_j \cdot \frac{(-1)^{j+i+1}a_{ij} \det(A_{ \hat{i}, \hat{j}})}{\det(A)},1\leq i\leq l \end{equation}
where $\det(A_{ \hat{i}, \hat{j}})$ denotes the $(l-1)\times(l-1)$ minor obtained by deleting $i$-th row and $j$-th column from $A$. Substituting (\ref{daititi}), (\ref{cramersolution}) becomes

\begin{equation}\label{kuaihaole} {x_{I_{ij}}} = {x_I}  \frac{(-1)^{j+i+1}a^3_{ij} \det(A_{ \hat{i}, \hat{j}})}{\det(A)},1\leq i\leq l.\end{equation}

Now to prove (\ref{zhenshia}), it suffices to show that $\frac{c_{I_{ij}}}{c_I}$ equals the same scalar factor as in the RHS of (\ref{kuaihaole}).
The expression of $c_I=c_{(1,2,\cdots,l)}$ is simple which is just $(-1)^{1+\cdots+l+1}\det(A)\det^3(I_l)\linebreak
=(-1)^{1+2+\cdots+l+1}\det(A)$. For $c_{I_{ij}}=c_{(1,2,\cdots,\hat{i},\cdots,l,l+j)}=(-1)^{1+2+\cdots+l-i+l+j+1}D_{I_{ij}}D^3_{ {I^c_{ij}}}$,   $D_{I_{ij}}$ is an $l\times l$ minor whose last row is $(0,\cdots,1\cdots,0)$ with $1$ lying on the $j$-th place, hence Laplace expansion along this row  gives
$D_{I_{ij}}=(-1)^{l+j}\det(A_{ \hat{i}, \hat{j}}).$
The minor $D_{ {I^c_{ij}}}$ has $(a_{i1},\cdots,a_{il})$ as the first row and the entry $a_{ij}$  therein is the only nonzero entry on the $j$-th column, thus applying Laplace expansion to $D_{ {I^c_{ij}}}$ along the $j$-th column gives
 $D_{{I^c_{ij}}}=(-1)^{j+1}a_{ij}.$

Finally $$\frac{c_{I_{ij}}}{c_I}=\frac{(-1)^{1+2+\cdots+l-i+l+j+1}(-1)^{l+j}\det(A_{ \hat{i}, \hat{j}})(-1)^{3{(j+1)}}a^3_{ij}}{(-1)^{1+2+\cdots+l+1}\det(A)}$$ which, after clearing powers of $-1$,  exactly equals  $\frac{(-1)^{j+i+1}a^3_{ij} \det(A_{ \hat{i}, \hat{j}})}{\det(A)}$.

\bigskip


For a general $M$ written as      $M=\scriptsize\left(
                              \begin{array}{c}
                                A \\
                                B \\
                              \end{array}
                            \right)$ \normalsize  where $A$ and $B$ are upper and lower $l\times l$ blocks,  then both $A$ and $B$ are invertible since $M$ is in general row position.  Let $\widetilde{M}=\left( \begin{smallmatrix} AB^{-1} \\ I \end{smallmatrix} \right)$ and $\{\widetilde{g}_i\}_{i=1}^{2l}$ be linear forms corresponding to rows of $\widetilde{M}$. Since $M=\widetilde{M}B$, it holds that $g_I(\mathbf{z})=\widetilde{g}_I(B\mathbf{z})$ for all $\mathbf{z}$ hence  by (\ref{0he}),  $\sum_I x_I\widetilde{g}_I^2(B\mathbf{z})=0.$ As $B$ is invertible, this  in turn implies  $\sum_I x_I\widetilde{g}^2_I=0.$

By what we have just proved in the standard block case, when $I=(1,2,\cdots,l)$ and $I_{ij}$ are as   (\ref{iij}), it holds that
\begin{equation}\label{zhenshiaaa}\frac{x_I}{\widetilde{c}_I}=\frac{x_{I_{ij}}}{\widetilde{c}_{I_{ij}}}.
\end{equation}
 where $\widetilde{c}_I,\widetilde{c}_{I_{ij}}$ are the corresponding \pkk polynomial coefficients with respect to $\widetilde{M}$. Finally, the desired (\ref{zhenshia}) follows from (\ref{zhenshiaaa}) and the   fact that(see Lemma  \ref{huanba} below) $$c_I={\det}^4(B)\widetilde{c}_I,  {c}_{I_{ij}}={\det}^4(B)\widetilde{c}_{I_{ij}}.$$
\end{proof}

\section{Structure  theory     of  \pkk polynomials}
 In this section we start  a  general structure theory  on $P_M$.
 In Section 3.1 we present elementary   facts  on how linear operation on $M$ affects the behaviour of $P_M$. In Section 3.2 we   prove a basic recursive formula which represents a $2l\times l$ standard \pkk polynomial by $2k\times k$ ones ($k<l$) along coordinate planes, and this implements an analogue of Laplace determinant expansion formula.      Section 3.3 consists of an example where we combine materials in Section 3.1 and 3.2 to  visualize how $P_M$ deviates from the zero polynomial as $l$ goes up via its full factorizability into linear forms.

    \subsection{Elementary properties}

 \begin{lem}\label{hanglieshi}
Let $M$ and $P_M$ be as in Definition \ref{pksum}.

 \begin{itemize}

 \item[(i)]
   Multiplying a  row of $M$ by a scalar $\lambda$ results in $P_M$ being scaled by $\lambda^3$;

 \item[(ii)] Swapping two rows of $M$ changes the sign of $P_M$.

\end{itemize}
Consequently,  $P_M$ vanishes if $M$ admits two proportional rows.
 \end{lem}

\begin{proof}
 $(i)$ is straightforward from Definition \ref{pksum} so we only prove $(ii)$. Since any row swap is the consequence of an odd number of adjacent row swaps, it suffices to prove $(ii)$ for adjacent row swaps. In the following proof we adopt the notation that for any $I=(i_1,\cdots,i_l)$,   $g_I:=g_{i_1} g_{i_2} \cdots g_{i_l} $ and $M_I$ denotes the $l\times l$ submatrix of $M$ using rows from $I$.

Now we fix $i\in\{1,2,\cdots,2l-1\}$.   Let $\widetilde{M}$ be
the matrix obtained by swapping  rows $i,i+1$ of $M$ and let $\sigma$ denote the permutation of $\{1,2,\cdots, 2l\}$  induced by the swap $i\leftrightarrow i+1$.


Writing  $$P_{\widetilde{M}}:=\sum_{I=(i_1,i_2\cdots, i_l)}(-1)^{|I|+1}\widetilde{D}_I\widetilde{D}^3_{I^c}\widetilde{g}_{I}^2,$$ if  we can  prove  that for any $I=(i_1,i_2\cdots, i_l) $,
\begin{equation}\label{swaap}(-1)^{|I|+1}D_ID^3_{I^c}g_{I}^2 =-((-1)^{|\sigma I|+1}\widetilde{D}_{\sigma I}\widetilde{D}^3_{({\sigma I})^c}\widetilde{g}_{\sigma I}^2 )\end{equation}
then taking sum over all $I$ gives $(ii)$  since $\sigma $ is a bijection.

We explicitly prove (\ref{swaap}) in the following cases which exhaust  all indices.

If $i\in I,i+1\in I$, then $I=\sigma I$, $M_{I^c}=\widetilde{M}_{({\sigma I})^c}$, $g_I=\widetilde{g}_{\sigma I}$, while $M_I$ differs from $\widetilde{M}_{\sigma I}$ by a row swap. Hence $|I|=|\sigma I|$, $D_{I^c}=\widetilde{D}_{({\sigma I})^c}$ and $D_I=-\widetilde{D}_{\sigma I}$, which implies (\ref{swaap}).

If $i\notin I,i+1\notin I$, then $I=\sigma I$, $M_I=\widetilde{M}_{\sigma I}$,  $g_I=\widetilde{g}_{\sigma I}$, while $M_{I^c}$ differs from $\widetilde{M}_{({\sigma I})^c}$ by a row swap. Hence $|I|=|\sigma I|$, $D_{I }=\widetilde{D}_{({\sigma I}) }$ and $D_{I^c}=-\widetilde{D}_{({\sigma I})^c}$, which also gives  (\ref{swaap}).

If $i\in I,i+1\notin I$(or $i\notin I,i+1\in I$), then $M_I=\widetilde{M}_{\sigma I}$,$M_{I^c}=\widetilde{M}_{({\sigma I})^c}$,$g_I=\widetilde{g}_{\sigma I}$, hence
 $D_{I }=\widetilde{D}_{({\sigma I}) }$, $D_{I^c}=\widetilde{D}_{({\sigma I})^c}$. In this case, (\ref{swaap}) follows from  $|I|=|\sigma I|-1$(or $|I|=|\sigma I|+1$).
\end{proof}

%
With respect to matrix multiplication, we have  the following  ``change of coordinates" formula     which follows immediately from    the determinant  product rule $det(AB)=det(A)det(B)$.
\begin{lem}\label{huanba}
Let $N$ be an $l\times l$ matrix and  let $M$ be a $2l\times l$ matrix. Then
  \begin{equation}\label{pmn}P_{MN}(\mathbf{z})=(\det N)^4P_M(N\textbf{z})\end{equation}
  where $N\mathbf{z}$ denotes the action of $N$ on the column vector $\mathbf{z}=(z_1,\cdots,z_l)^T$.\end{lem}


If $M$ can be blocked as $M=\scriptsize\left(
                              \begin{array}{c}
                                A \\
                                B \\
                              \end{array}
                            \right)$ \normalsize  such that  $B$ is   invertible, then applying (\ref{pmn})  gives

\begin{equation}\label{huanyuan}
P_M(\mathbf{z})  = (\det B)^4 P_{\tiny\left( \begin{smallmatrix} AB^{-1} \\ I \end{smallmatrix} \right)\normalsize} (B\mathbf{z}).
\end{equation}

This formula  will be used to deduce results on general \pkk polynomials from those on  standard ones.
We end this subsection with such a   lemma  on factorizability of \pkk polynomials to be used later.

\begin{lem}\label{zhege}Let $l\geq 2$ be a fixed integer. The following  are equivalent.

$(i)$  For any  $2l\times l$   matrix $M$ of the block form \scriptsize$\left(
                              \begin{array}{c}
                                A \\
                               I \\
                              \end{array}
                            \right)$\normalsize, $P_M$ admits $z_1z_2\cdots z_l$ as a factor.

$(ii)$ For any $2l\times l$   matrix $M$, $P_M$ admits $\prod_{i=1}^{2l}g_i$ as a factor, where $\{g_i\}_{i=1}^{2l}$ are linear forms corresponding to rows of $M$.

\end{lem}

 \begin{proof}It suffices to prove $(i)\Rightarrow(ii)$ as the converse direction is trivial. Let $M$ and $\{g_i\}_{i=1}^{2l}$ be as in $(ii)$.
  By Lemma \ref{hanglieshi}, we may assume $g_i$ is not proportional to $g_j$ for any $i\neq j$,  hence it suffices to show $g_i|P_M$ for each $1\leq i\leq 2l$. By symmetry, we only   prove $g_{2l}|P_M$.

In $M$ there always exists   a non-singular $l\times l$ submatrix   containing its $2l$-th row. Otherwise every $D_ID^3_{I^c}=0$  and thus $P_M=0$.  After a suitable   row permutation (which does not affect factorizability),  we may assume  $M=\scriptsize\left(
                              \begin{array}{c}
                                A \\
                                B \\
                              \end{array}
                            \right)$ \normalsize  such that  $B$ is   invertible and the last row of $B$ is  $g_{2l}$. Now    $z_l|P_{\left( \begin{smallmatrix} AB^{-1} \\ I \end{smallmatrix} \right)}$ by  $(i)$,  which   combined with   (\ref{huanyuan})  gives the desired  $ g_{2l}|P_M$  since $g_{2l}$ is exactly   $z_l$ composed with the action of $B$.
 \end{proof}

\subsection{Expansion  formula }

   In classical \pkk theory, $2l\times l$ matrices of the block form \scriptsize $\left(\begin{matrix} A \\ I_l \end{matrix}\right)$ \normalsize  are of particular interest where the \pkk relation boils down to Laplace expansion of $\det(A)$.
Our fundamental structural theorem on $P_M$ is the following analogue of the  Laplace determinant expansion formula, which says    that for any $k<l$,  a   $2l\times l $ standard  \pkk polynomial can   be    expanded   by        $2k\times k$  pieces along   coordinate $k$-planes.
\begin{thm}\label{lap}Let $l,m,k$ be positive integers such that $l=m+k$. Let  $A$ be an $l\times l$ matrix,  $\{g_i\}_{i=1}^l$ be linear forms corresponding to   rows of $A$, and   $M$ be a $2l\times l$ matrix of the block form $M=\scriptsize\left(
                              \begin{array}{c}
                                A \\
                                I_l \\
                              \end{array}
                            \right)$\normalsize.  It holds that

\begin{equation}\label{lape}
P_M(z_1,\cdots,z_k,0,\cdots,0)= \sum_{\tiny K=(k_1,k_2,\cdots,k_m)\subseteq (1,2,\cdots,l)\normalsize} (-1)^{|K|}D_K ({g}_{k_1}^2 {g}_{k_2}^2\cdots  {g}_{k_m}^2)|_{(z_1, \cdots,z_k,0,\cdots,0)}  P_{M_K}(z_1,\cdots,z_k).\end{equation}
Here  for a fixed $m$-index $K=(k_1,k_2,\cdots,k_m)$,   $|K|=k_1+\cdots+k_m$; $D_K$ denotes the determinant of the $m\times m$ matrix living in  $k_1,k_2\cdots, k_m$ rows  and the last $m$ columns of $A$; $M_K$ denotes the $2k\times k$ matrix $ \scriptsize\left(
                              \begin{array}{c}
                                A_K \\
                                I_{k } \\
                              \end{array}
                            \right)$\normalsize  where $A_K$ is the $k\times k$ submatrix of $A$    complementary to $D_K$.
\end{thm}

\bigskip

%

 Let $l,k,m$ be as in Theorem \ref{lap}.  We fix some notations   to be used in the forthcoming proof of  Theorem \ref{lap}.
 \begin{itemize}
 \item  $M_k$:   the $2l\times k$  matrix   consisting of the first $k$ columns of $M$ .

 \item $M_m$:  the $2l\times m$   matrix consisting of the last $m$ columns of $M$.
\item $ \Delta_K$:  the set of all $l$-indices $ I=(i_1,i_2\cdots, i_l)$ such that  $K\subseteq I\subseteq (1,\cdots,l+k)$, where $ K=(k_1,k_2,\cdots, k_m)\subseteq\{1,2,\cdots,l\}$ is  a fixed $m$-index.

\end{itemize}
The following notations make sense for   $I\in \Delta_K$.
\begin{itemize}
\item
 $M_{I_K}$:  the $k\times k$ matrix using  rows $I\setminus K$  and first $k$ columns of $M$(which is a submatrix of $M_K$).

   \item $M_{I^c_K}$:  the $k\times k$ matrix which is  complementary to $M_{I_K}$ in $M_K$.

  \item $D_{I_K}$ and $D_{I^c_K}$: the determinants of
   $M_{I_K}$ and $M_{I^c_K}$  respectively.

  \item $\sigma(I)$:  the $k$-index out of $\{1,2,\cdots,2k\}$   whose entries are row indices of   $M_{I_K}$ with respect to   $M_K$.
\end{itemize}

The following  example helps  to visualize above notations.
\begin{exa} Assume  $l=4,m=2,k=2$ and write  $M$  as
 \scriptsize $\begin{pmatrix}
a_1 & a_2 & a_3 & a_4 \\
b_1 & b_2 & b_3 & b_4 \\
c_1 & c_2 & c_3 & c_4 \\
d_1 & d_2 & d_3 & d_4 \\
1   & 0   & 0   & 0   \\
0   & 1   & 0   & 0   \\
0   & 0   & 1   & 0   \\
0   & 0   & 0   & 1
\end{pmatrix}.$ \normalsize

Take $$K=(2,3),$$  then $D_K=\scriptsize\left|
    \begin{array}{cc}
      b_3 &b_4 \\
      c_3 & c_4 \\
    \end{array}
  \right|$ \normalsize, $A_K=\scriptsize\left(
    \begin{array}{cc}
      a_1 & a_2 \\
      d_1 & d_2 \\
    \end{array}
  \right)$ \normalsize  and $M_K=\scriptsize \begin{pmatrix}
a_1 & a_2 \\
d_1 & d_2 \\
1   & 0   \\
0   & 1
\end{pmatrix}.$  \normalsize Take $I=(1,2,3,6)$, then  $I\in\Delta_K$ and  $M_{I_K}=\scriptsize\left(
    \begin{array}{cc}
      a_1 &a_2 \\
      0 & 1\\
    \end{array}
  \right)$ \normalsize,  $M_{I^c_K}=\scriptsize\left(
    \begin{array}{cc}
      d_1 &d_2 \\
      1 & 0\\
    \end{array}
  \right)$.  Since $M_{I_K}$ occupies the first and fourth rows of $M_K$, $\sigma(I)=(1,4)$.

 \end{exa}


The following    facts   will   be used in the proof of Theorem \ref{lap}.

 \begin{itemize}

 \item[(a)]    When $I$ exhausts $\Delta_K$ which has cardinality   $\binom{l+k-m}{l-m}= \binom{2k}{k}$,  $(M_{I_K},M_{I^c_K})$ exhausts all possible partitions of $M_{K}$ into two $k\times k$ matrices, and $\sigma(I)$ exhausts all   $k$-indices from $\{1,2,\cdots,2k\}$.

  \item[(b)]  For each
$I\subseteq (1,\cdots,l+k)$, the determinant $D_I$ can be written as  a linear sum of $\{D_K,   K \subseteq (1,2,\cdots,l)\}$ via Laplace expansion.

In fact,  applying Laplace expansion of   $D_I$
  along its last $m$ columns expresses  it  as a signed sum  $\sum \pm det(X)det(Y)$ where $X$ is  an $m\times m$ minor from $M_m$  and $Y$ is a $k\times k$ minor from $M_k$. As the lower $l\times l$ block of  $M$ is the identity matrix, the $k\times m$ matrix using   $l+1,\cdots,l+k$ rows and the last $m$ columns of $M$ is the zero matrix. Since $I\subseteq (1,\cdots,l+k)$,   $X$ will have a zero row  if it overlaps the lower $l\times l$ block of $M$. In other words,   $det(X)$ can be nonzero only  when $X$   lies entirely in the upper $l\times l$ block of $M$, hence  $det(X)=D_K$ for some $K=
(k_1,k_2\cdots, k_m)\subseteq(1,2,\cdots,l)$.

  \item[(c)] For each $I\in\Delta_K$, $D_{I^c_K}=D_{I^c}.$

  This follows from the fact that the $l\times l$ matrix corresponding to $I^c$ is of the block form $\left(
    \begin{array}{cc}
      M_{I^c_K} &* \\
      0 & I_{m\times m}\\
    \end{array}
  \right).$
 \end{itemize}

\begin{proof}


Since  $ g_{l+1}=z_1, g_{l+2}=z_2,\cdots,g_{2l}=z_l,$   we have $g_{l+k+1}=\cdots=g_{2l}=0$   restricted to $(z_1, \cdots,z_k, 0,\cdots,0)$, so the $I$-summand \begin{equation}\label{snd}(-1)^{|I|+1}D_ID^3_{I^c}g_{i_1}^2g_{i_2}^2\cdots g_{i_l}^2\end{equation}  in   (\ref{pmii}) vanishes if $I$ has entries coming from $\{ l+k+1,l+k+2,\cdots,2l\}$, which implies
\begin{equation}\label{nbv}P_M(z_1,z_2,\cdots,z_k,0,\cdots,0)=\sum_{I=(i_1,i_2\cdots, i_l)\subseteq\{1,\cdots,l+k\}}(-1)^{|I|+1}D_ID^3_{I^c}g_{i_1}^2g_{i_2}^2\cdots g_{i_l}^2|_{(z_1, \cdots,z_k,0,\cdots,0)}\end{equation}

By above fact $(b)$, applying Laplace expansion to each $D_I, I \subseteq\{1,\cdots,l+k\}$ along the last $m$ columns expresses the entire $I$-summand in (\ref{nbv})  as a ``linear sum" of $\{D_K,   K \subseteq\{1,2,\cdots,l\}\}$ hence the entire $P_M(z_1,z_2,\cdots,z_k,0,\cdots,0)$ can be expressed as
 \begin{equation}\label{biaoshiling}P_M(z_1,z_2,\cdots,z_k,0,\cdots,0)=\sum_{K=
(k_1,k_2\cdots, k_m)\subseteq(1,2,\cdots,l)} D_K\Lambda_K \end{equation} where $\Lambda_K$ is the  ``polynomial coefficient" of $D_K$.
%

Now it  suffices to show that for any fixed    $m$-index  $K=(k_1,k_2,\cdots,k_m)\subseteq(1,2,\cdots,l) $,   $\Lambda_K$    can be expressed as  \begin{equation}\label{spkm}\Lambda_K=(-1)^{|K|}  ({g}_{k_1}^2 {g}_{k_2}^2\cdots  {g}_{k_m}^2)|_{(z_1, \cdots,z_k,0,\cdots,0)}  P_{M_K}(z_1,\cdots,z_k).\end{equation}  To this end,  we need  to  trace the $D_K$-component  in each $I$-summand of (\ref{nbv}) and sum these components together.
More precisely, $D_K$ appears in the Laplace expansion of $D_I$ only when    $K\subseteq I$ hence it suffices to trace over $I\in \Delta_K$.

 Since  $M_{I_K}$ is   exactly the complementary matrix of $D_K$ in $D_I$, the $D_K$ contribution to $D_I$ is the Laplace expansion summand $\epsilon(I,K) D_KD_{I_K}$, where $\epsilon(I,K)=\pm 1$ is the sign in the expansion.
Inserting this and $D_{I^c_K}=D_{I^c}$(by above fact $(c)$),  we see that the  $D_K$-component  of  the $I$-summand in (\ref{nbv})  is the following expression


\begin{equation}\label{ssnd}(-1)^{|I|+1}\epsilon(I,K) (D_KD_{I_K})D^3_{I^c_K}g_{i_1}^2g_{i_2}^2\cdots g_{i_l}^2=D_Kg_{K}^2[ (-1)^{|I|+1} \epsilon(I,K) D_{I_K}D^3_{I^c_K}g_{I\setminus K}^2]\end{equation} restricted to $(z_1, \cdots,z_k,0,\cdots,0)$,
where $g_{K}^2$, $ g_{I\setminus K}^2$ denote   the obvious split of  $ g_{i_1}^2g_{i_2}^2\cdots g_{i_l}^2$ according to the partition    $K\cup (I\setminus K)$ of $I$.

Taking sum over all $I\in \Delta_K$ and factoring out $D_K$ gives

\begin{equation}\label{fuhao}\Lambda_K= g_{K}^2|_{(z_1, \cdots,z_k,0,\cdots,0)} \sum_{I\in \Delta_K} (-1)^{|I|+1} \epsilon(I,K) D_{I_K}D^3_{I^c_K}g_{I\setminus K}^2|_{(z_1, \cdots,z_k,0,\cdots,0)}. \end{equation}

On the other hand,
  $M_{I_K}$ uses rows $I\setminus K$  and first $k$ columns of $M$, hence  linear forms in  $g_{I\setminus K}^2|_{(z_1, \cdots,z_k,0,\cdots,0)}$   exactly correspond to rows of $M_{I_K}$, hence we have, by Definition \ref{pksum} and above fact $(a)$,
  \begin{equation}\label{mkdef} P_{M_K}(z_1,\cdots,z_k)=\sum_{I\in \Delta_K} (-1)^{|\sigma(I)|+1} D_{I_K}D^3_{I^c_K}g_{I\setminus K}^2|_{(z_1, \cdots,z_k,0,\cdots,0)}.\end{equation}

  Comparing (\ref{mkdef}) and (\ref{fuhao}), we see that the desired   (\ref{spkm}) holds if the following sign  relation holds
  \begin{equation}\label{match}
  (-1)^{|I|+1} \epsilon(I,K)= (-1)^{|\sigma(I)|+1} (-1)^{|K|},\forall I\in\Delta_K.
  \end{equation}


\bigskip

The remaining proof consists of two steps.
We first verify (\ref{match}) assuming  $K=(1,2,\cdots, m)$, which will imply (\ref{spkm}) for this particular $K$.  For general $K$, we will not directly check (\ref{match}) which, though  not hard, is  quite tricky. Instead,  we prove (\ref{spkm}) by reducing to $K=(1,2,\cdots, m)$ via row permutation.


  \bigskip

 Now we begin with $K=(1,2,\cdots, m)$.  For this $K$, any $I\in\Delta_K$ can be written as $(1,2,\cdots,m,t_1,\cdots,t_k)$ where $ t_i>m,i=1,2,\cdots,k$,  thus $\sigma(I)=(t_1-m,t_2-m,\cdots,t_k-m)$.  Since $D_K$ occupies rows $1,2,\cdots,m$ and columns  $(k+1),(k+2),\cdots,(k+m)$ of $D_I$, the Laplace sign $\epsilon(I,K)$ is $$(-1)^{(k+1)+(k+2)+\cdots+(k+m)+1+2+\cdots+m}=(-1)^{mk}.$$
 Now the LHS of (\ref{match}) equals
 $$(-1)^{mk+(1+2+\cdots+m+t_1+\cdots+t_k)+1}$$ while the RHS equals
 $$(-1)^{ 1+(t_1-m)+(t_2-m)+\cdots+(t_k-m) +1+2+\cdots+m},$$
 thus $\frac{\mathrm{LHS}}{\mathrm{RHS}}=(-1)^{2mk}=1$ as desired.

 \bigskip

Finally   we prove  (\ref{spkm}) for an arbitrary fixed $K=(k_1,k_2,\cdots,k_m)$, whose  complementary index  in $(1,2,\cdots,l)$ is denoted by  $(t_1,t_2,\cdots, t_k)$.


Let  $\widetilde{M}=\scriptsize\left(
                              \begin{array}{c}
                                \widetilde{A} \\
                                I \\
                              \end{array}
                            \right)$ \normalsize,   where  $\widetilde{A}$ denotes  the $l\times l$ matrix  obtained by permuting rows of the original matrix $A$   so that  rows $ k_1,k_2,\cdots,k_m $ of $A$ occupy  rows $1,2,\cdots,m$  of $\widetilde{A}$, and  rows  $ t_1,t_2,\cdots, t_k $ of $A$ occupy      rows $m+1,\cdots,l$ of $\widetilde{A}$.

                              Now    it holds that $D_K=\widetilde{D}_{(1,2,\cdots,m)}, M_K=\widetilde{M}_{(1,2,\cdots,m)}$ and
                            $$g_{k_i}=\widetilde{g}_i,1\leq i\leq m.$$
                             Let  $\epsilon(A,\widetilde{A})$ denote the sign of this row   permutation,
then by Lemma \ref{hanglieshi},
$$P_M=\epsilon(A,\widetilde{A})P_{\widetilde{M}}.$$
If we write \begin{equation} P_{\widetilde{M}}(z_1,z_2,\cdots,z_k,0,\cdots,0)=\sum_{K=
(k_1,k_2\cdots, k_m)\subseteq(1,2,\cdots,l)} \widetilde{D}_K\widetilde{\Lambda}_K,\end{equation}
then  $\Lambda_K=\epsilon(A,\widetilde{A})\widetilde{\Lambda}_{(1,2,\cdots,m)}$.


 By validity of (\ref{spkm}) for $(1,2,\cdots,m)$,
 \begin{align}
\text{\normalsize $\displaystyle
\begin{aligned}
    \widetilde{ \Lambda}_{(1,2,\cdots,m)}&=(-1)^{| (1,2,\cdots,m)|}  ({\widetilde{g}}_{ 1}^2 {\widetilde{g}}_{ 2}^2\cdots  {\widetilde{g}}_{ m}^2)|_{(z_1, \cdots,z_k,0,\cdots,0)}  P_{\widetilde{M}_{(1,2,\cdots,m)}}(z_1,\cdots,z_k) \\
    &=(-1)^{1+2+\cdots+m} ({g}_{k_1}^2 {g}_{k_2}^2\cdots  {g}_{k_m}^2)|_{(z_1, \cdots,z_k,0,\cdots,0)}  P_{M_K}(z_1,\cdots,z_k)
\end{aligned}
$}
\end{align}
thus
 $$\Lambda_K=\epsilon(A,\widetilde{A})(-1)^{1+2+\cdots+m}  ({g}_{k_1}^2 {g}_{k_2}^2\cdots  {g}_{k_m}^2)|_{(z_1, \cdots,z_k,0,\cdots,0)}  P_{M_K}(z_1,\cdots,z_k).  $$  Finally,  the desired
(\ref{spkm})  holds for $K=(k_1,k_2,\cdots,k_m)$  provided that $$\epsilon(A,\widetilde{A})(-1)^{1+2+\cdots+m}=(-1)^{|K|},  $$ and this indeed holds since
 there are exactly    $k_1+\cdots+k_m-(1+2+\cdots+m)$ inversions  in the permutation $\{k_1,\cdots,k_m,t_1,\cdots,t_k\} $ of $\{1,2,\cdots,l\}$.


\end{proof}

The following corollary to Theorem \ref{lap}   will be used in Section 3.3.
 \begin{cor}\label{zhengchu}
Let   $l\geq 2$ be a fixed integer. If all $2(l-1)\times (l-1)$  {Pl$\ddot{ {u}}$cker   polynomials} are identically zero, then for any $2l\times l$ matrix $M$, $P_M$ admits $\prod_{i=1}^{2l}g_i$ as a factor, where $\{g_i\}_{i=1}^{2l}$ are row linear forms of $M$.

\end{cor}
\begin{proof}Since all $2(l-1)\times (l-1)$  {Pl$\ddot{ {u}}$cker   polynomials} are identically zero, Theorem \ref{lap} implies that $P_M|_{z_i=0}=0, 1\leq i\leq l$ for any  standard $2l\times l$ \pkk polynomial $P_M$, and the conclusion follows from Lemma \ref{zhege}.

\end{proof}


 \subsection{Application to factorizability}

Since  $P_M$ is a high-degree polynomial, it cannot be completely determined by its value on coordinate planes, thus formula (\ref{lape}) alone is not sufficient to fully capture     a $2l\times l$ \pkk polynomial even if everything on $2(l-1)\times (l-1)$ ones is known.
In fact, the structural complexity of $P_M$ indeed increases very rapidly as $l$ goes up, and we present  an example   to   visualize   this across  $l=2,3,4$, which  exhibits  an elegant inductive picture  on full factorizability of $P_M$ into linear forms.

\begin{exa}\label{36}  Let $M=\{a_{ij}\}_{1\leq i\leq 2l, 1\leq j\leq l}$  be a $2l\times l$ matrix and $g_i= a_{i1}z_1+\cdots+a_{il}z_l,  1\leq i\leq 2l$. Then

$(i)$ When $l=2$, $P_M=0.$

$(ii)$ When $l=3$,  $P_M$ admits $\prod_{i=1}^{6}g_i$ as a factor, and   $P_M=0$ if  $\{g_i^2\}^{6}_{i=1}$ are linearly dependent.

$(iii)$ When $l=4$,  $P_M$ admits $\prod_{i=1}^{8}g_i$ as a factor if  $\{g_i^2\}^{8}_{i=1}$ are linearly dependent.
    \end{exa}

 \begin{proof}

$(i)$ Write $M=\tiny\left(
\begin{array}{cccc}
 a_1 & b_1 & c_1 & d_1 \\
 a_2 & b_2 & c_2 & d_2 \\
\end{array}
\right)^T$\normalsize.
We first    show that  $P_M$ admits $(a_1z_1 + a_2z_2) (b_1z_1 + b_2z_2)(c_1z_1 + c_2z_2) (d_1z_1 + d_2z_2)$ as a factor. Since $P_M$ is of degree 4 in $\{z_1,z_2\}$, this will imply    \begin{equation}\label{22fac}P_M=Q(a_1z_1 + a_2z_2) (b_1z_1 + b_2z_2)(c_1z_1 + c_2z_2) (d_1z_1 + d_2z_2)\end{equation} where $Q$ is a polynomial involving only  $a_i,b_i,c_i,d_i$. Then we show $Q=0$ to finish.

The first part follows from the fact that $2\times 1$ \pkk polynomials are identically zero. Precisely, for any $2\times 1$ matrix $\scriptsize \left(
                              \begin{array}{c}
                                a  \\
                                b \\
                              \end{array}
                            \right)$\normalsize, $P_{\tiny\left(
                              \begin{array}{c}
                                a  \\
                                b  \\
                              \end{array}
                            \right)}=ab^3(az_1)^2-a^3b(bz_1)^2$ which is trivially zero. Now  by Corollary \ref{zhengchu}, $(a_1z_1 + a_2z_2) (b_1z_1 + b_2z_2)(c_1z_1 + c_2z_2) (d_1z_1 + d_2z_2)|P_M$.

     It remains to show that the factor $Q$ in (\ref{22fac}) is the zero polynomial. Recall that $P_M$ is of total degree 3 in each row of $M$ (by Lemma  \ref{hanglieshi}) while one degree is absorbed  by $(a_1z_1 + a_2z_2) (b_1z_1 + b_2z_2)(c_1z_1 + c_2z_2) (d_1z_1 + d_2z_2) $, hence $Q$ is of degree 2 in each row of $M$.
     For convenience we use the same notation $\{g_i\}_{i=1}^4$ to denote the 4 rows of $M$.

      Let $\{q_i\}_{i=1}^3$ denote the   standard  basis   spanning the space of quadratic forms in two variables, that is,
\[
q_1=x_1^2,\quad q_2=x_2^2,\quad
q_3=x_1x_2,
\]
then   $q_j(g_i)$ makes sense as a monomial in the two entries of $g_i$. As $Q$ is quadratic  in each $g_i$, expanding   with respect to the basis  $\{q_j\}_{j=1}^3$ gives the following  expression of  $Q$    as a   sum over all $3^4$ indices $(r_1,r_2,r_3,r_4)$,  $1\leq r_i\leq 3$
\begin{equation}\label{44sum}
Q(g_1,g_2,g_3,g_4)
=
 \sum_{\stackrel{(r_1,r_2,r_3,r_4)}{\quad }} c_{r_1,r_2,r_3,r_4}\,
q_{r_1}(g_1)q_{r_2}(g_2)q_{r_3}(g_3)q_{r_4}(g_4),
\end{equation}where $c_{r_1,r_2,r_3,r_4}$  are scalar coefficients.
Since $P_M$ is anti-symmetric in $(g_1,g_2,g_3,g_4)$  while its factor $\prod_{i=1}^4 g_i$ is symmetric,    $Q$ as the remaining factor has to be anti-symmetric,  which in turn gives
$$c_{ \dots r_{i } \dots r_{j }\dots }=-c_{ \dots r_{j } \dots r_{i }\dots } $$ for any $1\leq i <j \leq 4$. In particular, $c_{r_1,r_2,r_3,r_4}=0$   if     $r_i=r_j$ for some $i\neq j$, and this is inevitable by pigeon hole   as each $r_i$ lies in the set $\{1,2,3\}$.\bx

\bigskip



%
$(ii)$
 Write
$ M=
\tiny \left(
\begin{array}{cccccc}
 a_1 & b_1 & c_1 & d_1 & e_1 & f_1 \\
 a_2 & b_2 & c_2 & d_2 & e_2 & f_2  \\
 a_3 & b_3 & c_3 &  d_3 & e_3 & f_3 \\
\end{array}
\right)^T.$ We show that the following algebraic identity holds, from which   the $l=3$ part of Theorem \ref{36} follows.\begin{equation}\label{niubiplus}
P_M=-6 \mbox{\scriptsize $ \left| \begin{array}{cccccc}
a_1^2 & a_2^2 & a_3^2 & a_1 a_2 & a_1 a_3 & a_2 a_3 \\
b_1^2 & b_2^2 & b_3^2 & b_1 b_2 & b_1 b_3 & b_2 b_3 \\
c_1^2 & c_2^2 & c_3^2 & c_1 c_2 & c_1 c_3 & c_2 c_3 \\
d_1^2 & d_2^2 & d_3^2 & d_1 d_2 & d_1 d_3 & d_2 d_3 \\
e_1^2 & e_2^2 & e_3^2 & e_1 e_2 & e_1 e_3 & e_2 e_3 \\
f_1^2 & f_2^2 & f_3^2 & f_1 f_2 & f_1 f_3 & f_2 f_3
\end{array} \right|
(\sum_{i=1}^3a_iz_i)(\sum_{i=1}^3b_iz_i)(\sum_{i=1}^3c_iz_i)(\sum_{i=1}^3d_iz_i) (\sum_{i=1}^3e_iz_i)(\sum_{i=1}^3f_iz_i) $}
\end{equation}

In fact, given this identity,  $\prod^6_{i=1} g_i|P_M$ follows immediately. Moreover, $P_M$ vanishes if and only if  the  $6\times 6$ determinant in RHS of (\ref{niubiplus})  vanishes, and vanishing of this determinant is exactly the testimony for linear dependence of    $ \footnotesize\{(\sum_{i=1}^3a_iz_i)^2, \cdots,  (\sum_{i=1}^3f_iz_i)^2\}$ \normalsize.

Since  all $4\times 2$ \pkk polynomials are identically zero as we have already proved,  Corollary \ref{zhengchu} gives  $ \prod^6_{i=1} g_i|P_M$ where $\{g_i\}^6_{i=1}$ are linear forms corresponding to rows of $M$. Thus $P_M$ can be written as  $ P_M=Q\prod^6_{i=1} g_i$  where $Q$  is a  polynomial   involving  only $\{a_i,b_i,c_i,d_i,e_i,f_i,i=1,2,3\}$. A degree count forces $Q$  to be of degree 2 in $\{g_i\}_{i=1}^6$, and anti-symmetry of $P_M$ in $\{g_i\}_{i=1}^6$ forces $Q$ to be anti-symmetric as well.

  Let $\{q_i\}_{i=1}^6$ denote the   standard  basis   spanning the space of quadratic forms in three variables, that is,
\[
q_1=x_1^2,\quad q_2=x_2^2,\quad q_3=x_3^2,\quad
q_4=x_1x_2,\quad q_5=x_1x_3,\quad q_6=x_2x_3,
\]
then  $Q$   can be written  as a   sum over all $6^6$ indices $(r_1,r_2,r_3,r_4,r_5,r_6)$,  $1\leq r_i\leq 6$ as follows
\begin{equation}\label{66sum}
Q(g_1,\cdots,g_6)
=
 \sum_{\stackrel{(r_1,r_2,r_3,r_4,r_5,r_6)}{\quad }} c_{r_1\dots r_6}\,
q_{r_1}(g_1)q_{r_2}(g_2)q_{r_3}(g_3)q_{r_4}(g_4)q_{r_5}(g_5)q_{r_6}(g_6),
\end{equation}
where $c_{r_1\dots r_6}$  are scalar coefficients.
Set $\lambda:=c_{1,2,3,4,5,6}$. Anti-symmetry of $Q$ in  $(g_1,\cdots,g_6)$
%
%
gives \[
Q(g_1,\dots,g_6)
=
\lambda \sum_{\sigma\in S_6} \operatorname{sgn}(\sigma)\,
q_{\sigma(1)}(g_1)\cdots q_{\sigma(6)}(g_6)
=\lambda \scriptsize\left| \begin{array}{cccccc} a_1^2 & a_2^2 & a_3^2 & a_1 a_2 & a_1 a_3 & a_2 a_3 \\ b_1^2 & b_2^2 & b_3^2 & b_1 b_2 & b_1 b_3 & b_2 b_3 \\ c_1^2 & c_2^2 & c_3^2 & c_1 c_2 & c_1 c_3 & c_2 c_3 \\ d_1^2 & d_2^2 & d_3^2 & d_1 d_2 & d_1 d_3 & d_2 d_3 \\ e_1^2 & e_2^2 & e_3^2 & e_1 e_2 & e_1 e_3 & e_2 e_3 \\ f_1^2 & f_2^2 & f_3^2 & f_1 f_2 & f_1 f_3 & f_2 f_3 \end{array} \right|. \]\normalsize
That $\lambda=-6$ can be derived by assigning specific values, and we omit the elementary computation since this precise value is not critical for our conclusion.

\bigskip

Specializing (\ref{niubiplus}) to the standard case gives the following identity
  \begin{equation}\label{yibanban}\Large P_{\tiny\left(
                    \begin{array}{cccccc}
                      a_1 & b_1 & c_1 & 1 & 0 &0 \\
                      a_2 & b_2 & c_2 & 0 & 1 &0 \\
                      a_3 & b_3 & c_3 & 0 & 0 &1\\
                    \end{array}
                  \right)^T}\normalsize\text{\footnotesize  $\displaystyle =6\left|                  \begin{array}{ccc}
                    a_1a_2&b_1b_2 & c_1c_2\\
                    a_1a_3 & b_1b_3 & c_1c_3 \\
                    a_2a_3 & b_2b_3 & c_2c_3 \\
                  \end{array}
                \right|z_1 z_2 z_3 \left(a_1 z_1+a_2 z_2+a_3 z_3\right) \left(b_1 z_1+b_2 z_2+b_3 z_3\right) \left(c_1 z_1+c_2 z_2+c_3 z_3\right)$}.  \end{equation} We will use this identity as well as its $abd,acd,bcd$ versions in the following proof of the $l=4$ part.

\bigskip

$(iii)$
 In formula (\ref{huanyuan}),    linear forms corresponding to  rows of $M$ and  rows of  $\left( \begin{smallmatrix} AB^{-1} \\ I \end{smallmatrix} \right)$ are related by composition with the invertible map $B$, hence   both properties, $\prod_{i=1}^8 g_i|P_M$ and linear dependence of  $\{g_i^2\}^{8}_{i=1}$,     carry over from $M$  to $\left( \begin{smallmatrix} AB^{-1} \\ I \end{smallmatrix} \right)$ and vice versa. This implies that it suffices to prove the conclusion for standard $8\times 4$ \pkk polynomials. In light of Lemma \ref{zhege}, it suffices to  show that if

  \begin{equation}\label{8444}
M=\text{\scriptsize  $\displaystyle\left(
\begin{array}{cccccccc}
 a_1 & b_1 & c_1&d_1 & 1 & 0 & 0& 0 \\
 a_2 & b_2 & c_2&d_2 & 0 & 1 & 0 & 0\\
 a_3 & b_3 & c_3 &d_3& 0 & 0 & 1& 0 \\
 a_4 & b_4 & c_4 &d_4& 0 & 0 & 0& 1 \\
\end{array}
\right)^T$},  \end{equation} \normalsize  and
$$ \text{\footnotesize  $\displaystyle\{g_i^2\}_{i=1}^8= \{ (\sum_{i=1}^4a_iz_i)^2, (\sum_{i=1}^4b_iz_i)^2,
(\sum_{i=1}^4c_iz_i)^2,(\sum_{i=1}^4d_iz_i)^2,z_1^2
,z_2^2,z_3^2,z_4^2\}  $}$$   are linearly dependent, then $z_1z_2z_3z_4|P_M$.
Routine linear algebra shows  that above $\{g^2_i\}_{i=1}^8$ are linearly dependent if and only if
\begin{equation}\label{mat46}\operatorname{rank}W<4,\quad \mathrm{where}\quad W=\text{\footnotesize  $\displaystyle
 \left(
\text{\footnotesize $
\begin{array}{cccc}
 a_1 a_2 & b_1 b_2 & c_1 c_2 & d_1 d_2 \\
 a_1 a_3 & b_1 b_3 & c_1 c_3 & d_1 d_3 \\
 a_1 a_4 & b_1 b_4 & c_1 c_4 & d_1 d_4 \\
 a_2 a_3 & b_2 b_3 & c_2 c_3 & d_2 d_3 \\
 a_2 a_4 & b_2 b_4 & c_2 c_4 & d_2 d_4 \\
 a_3 a_4 & b_3 b_4 & c_3 c_4 & d_3 d_4
\end{array}
$} \right) .$}
\end{equation}
 In what follows we  prove that  if  $\operatorname{rank}W<4$, then $P_M|_{z_4=0}=0$,  and the conclusion follows by symmetry.

 By Theorem \ref{lap}, inserting $z_4=0$,   $P_M$ becomes  a linear form in $a_4,b_4,c_4,d_4$  whose  coefficients   as polynomials in other variables can respectively be divided by the LHS of   (\ref{yibanban})  and  its $bcd,acd,abd$ versions.
In particular, the $d_4$ coefficient is   $ (d_1 z_1+d_2 z_2+d_3 z_3) ^2$ times   $P_{\tiny\left(
                    \begin{array}{cccccc}
                      a_1 & b_1 & c_1 & 1 & 0 &0 \\
                      a_2 & b_2 & c_2 & 0 & 1 &0 \\
                      a_3 & b_3 & c_3 & 0 & 0 &1\\
                    \end{array}
                  \right)^T}$      hence by   (\ref{yibanban}),   this coefficient equals

\begin{equation}\label{d4}
\text{\footnotesize  $\displaystyle
 6 \begin{vmatrix}
    a_1a_2 & b_1b_2 & c_1c_2 \\
    a_1a_3 & b_1b_3 & c_1c_3 \\
    a_2a_3 & b_2b_3 & c_2c_3
\end{vmatrix} z_1 z_2 z_3 (a_1 z_1+a_2 z_2+a_3 z_3) (b_1 z_1+b_2 z_2+b_3 z_3) (c_1 z_1+c_2 z_2+c_3 z_3)(d_1 z_1+d_2 z_2+d_3 z_3)^2
$}
\end{equation}
In the same way, the coefficients of $a_4,b_4,c_4$    respectively equal
\begin{equation}\label{a4}
\text{\footnotesize $\displaystyle
-6 \begin{vmatrix}
    b_1b_2 & c_1c_2 & d_1d_2 \\
    b_1b_3 & c_1c_3 & d_1d_3 \\
    b_2b_3 & c_2c_3 & d_2d_3
\end{vmatrix} z_1 z_2 z_3 (a_1 z_1+a_2 z_2+a_3 z_3)^2 (b_1 z_1+b_2 z_2+b_3 z_3) (c_1 z_1+c_2 z_2+c_3 z_3)(d_1 z_1+d_2 z_2+d_3 z_3)
$}
\end{equation}

\begin{equation}\label{b4}
\text{\footnotesize  $\displaystyle
6 \begin{vmatrix}
    a_1a_2 & c_1c_2 & d_1d_2 \\
    a_1a_3 & c_1c_3 & d_1d_3 \\
    a_2a_3 & c_2c_3 & d_2d_3
\end{vmatrix} z_1 z_2 z_3 (a_1 z_1+a_2 z_2+a_3 z_3) (b_1 z_1+b_2 z_2+b_3 z_3)^2 (c_1 z_1+c_2 z_2+c_3 z_3)(d_1 z_1+d_2 z_2+d_3 z_3)
$}
\end{equation}

\begin{equation}\label{c4}
\text{\footnotesize  $\displaystyle
-6 \begin{vmatrix}
    a_1a_2 & b_1b_2 & d_1d_2 \\
    a_1a_3 & b_1b_3 & d_1d_3 \\
    a_2a_3 & b_2b_3 & d_2d_3
\end{vmatrix} z_1 z_2 z_3 (a_1 z_1+a_2 z_2+a_3 z_3) (b_1 z_1+b_2 z_2+b_3 z_3) (c_1 z_1+c_2 z_2+c_3 z_3)^2(d_1 z_1+d_2 z_2+d_3 z_3)
$}
\end{equation}
Now we have  \begin{equation}\label{p40}P_M|_{z_4=0} =d_4(\ref{d4})+c_4(\ref{c4})+b_4(\ref{b4})+a_4(\ref{a4}).\end{equation}\normalsize

After factoring out the common factor
\begin{equation}\label{cfac}6z_1 z_2 z_3\left(a_1 z_1+a_2 z_2+a_3 z_3\right) \left(b_1 z_1+b_2 z_2+b_3 z_3\right) \left(c_1 z_1+c_2 z_2+c_3 z_3\right)(d_1 z_1+d_2 z_2+d_3 z_3)\end{equation} \normalsize \normalsize of above four summands, proof of  $P_M|_{z_4=0}=0$  amounts to showing that


\begin{equation}\label{46}\text{\footnotesize  $\displaystyle
\begin{aligned}& d_4(d_1 z_1+d_2 z_2+d_3 z_3)\left|                  \begin{array}{ccc}
                    a_1a_2&b_1b_2 & c_1c_2\\
                    a_1a_3 & b_1b_3 & c_1c_3 \\
                    a_2a_3 & b_2b_3 & c_2c_3 \\
                  \end{array}
                \right|-c_4\left(c_1 z_1+c_2 z_2+c_3 z_3\right)\left|                  \begin{array}{ccc}
                    a_1a_2&b_1b_2 & d_1d_2\\
                    a_1a_3&b_1b_3 & d_1d_3 \\
                    a_2a_3&b_2b_3 & d_2d_3 \\
                  \end{array}
                \right|    \\&+b_4\left(b_1 z_1+b_2 z_2+b_3 z_3\right)\left|                  \begin{array}{ccc}
                    a_1a_2&c_1c_2 & d_1d_2\\
                    a_1a_3&c_1c_3 & d_1d_3 \\
                    a_2a_3&c_2c_3 & d_2d_3 \\
                  \end{array}
                \right|-  a_4\left(a_1 z_1+a_2 z_2+a_3 z_3\right) \left|                  \begin{array}{ccc}
                    b_1b_2&c_1c_2 & d_1d_2\\
                    b_1b_3&c_1c_3 & d_1d_3 \\
                    b_2b_3&c_2c_3 & d_2d_3 \\
                  \end{array}
                \right| \\& =0,  \end{aligned}$}\end{equation}
                or equivalently,

                \begin{equation}\label{nxz}\text{\footnotesize  $\displaystyle\left|
\begin{array}{cccc}
 a_1 a_2 & b_1 b_2 & c_1 c_2 & d_1 d_2 \\
 a_1 a_3 & b_1 b_3 & c_1 c_3 & d_1 d_3 \\
 a_2 a_3 & b_2 b_3 & c_2 c_3 & d_2 d_3 \\
 a_1 a_4 & b_1 b_4 & c_1 c_4 & d_1 d_4 \\
\end{array}
\right|z_1+\left|
\begin{array}{cccc}
 a_1 a_2 & b_1 b_2 & c_1 c_2 & d_1 d_2 \\
 a_1 a_3 & b_1 b_3 & c_1 c_3 & d_1 d_3 \\
 a_2 a_3 & b_2 b_3 & c_2 c_3 & d_2 d_3 \\
 a_2 a_4 & b_2 b_4 & c_2 c_4 & d_2 d_4 \\
\end{array}
\right|z_2+\left|
\begin{array}{cccc}
 a_1 a_2 & b_1 b_2 & c_1 c_2 & d_1 d_2 \\
 a_1 a_3 & b_1 b_3 & c_1 c_3 & d_1 d_3 \\
 a_2 a_3 & b_2 b_3 & c_2 c_3 & d_2 d_3 \\
 a_3 a_4 & b_3 b_4 & c_3 c_4 & d_3 d_4 \\
\end{array}
\right|z_3=0,$}\end{equation}\normalsize
but this  immediately follows from (\ref{mat46}) which forces   all $4\times 4$ minors of   the  matrix $W$  to be zero.  In (\ref{nxz})   three    minors are involved, and the vanishing of another 9 minors in (\ref{mat46})  will be  used for the proofs of $P_M|_{z_3=0}=0,$ $P_M|_{z_2=0}=0$ and $P_M|_{z_1=0}=0 $ in the same way. \end{proof}

\bigskip

SCHOOL OF MATHEMATICS, SHANDONG UNIVERSITY, JINAN 250100, CHINA.

\footnotesize \emph{Email:} \textbf{ lchencz@sdu.edu.cn}


\scriptsize

\end{document}